\documentclass[11pt]{article}

\usepackage{amssymb,amstext,amsthm,latexsym}
\usepackage
{amsmath}
\usepackage{amsfonts}

\usepackage{mathrsfs}
\usepackage[bbgreekl]{mathbbol}
\usepackage{url}
\usepackage{doi}
\usepackage{mparhack}

\theoremstyle{definition}

\theoremstyle{remark}

\numberwithin{equation}{section}

\input{64.sty}

\begin{document}

\title{A product forcing model in which the 
Russell-nontypical sets satisfy ZFC strictly 
between HOD and the universe\thanks
{Partial support of RFBR grant 20-01-00670 acknowledged.}}

\author{Vladimir~Kanovei\thanks 
{IITP RAS, Moscow, Russia, {\tt kanovei@iitp.ru}.}
\and
Vassily~Lyubetsky\thanks 
{IITP RAS, Moscow, Russia, {\tt lyubetsk@iitp.ru}}
}


\date{}

\maketitle

\begin{abstract}
A set is nontypical in the Russell sense, 
if it belongs to a countable ordinal definable set. 
The class $\HNT$ of all hereditarily nontypical sets 
satisfies all axioms of $\ZF$ and the double inclusion
$\HOD\sq\HNT\sq\rV$ holds. 
Solving a problem recently proposed by Tzouvaras, 
a generic extension $\rL[a,x]$ of $\rL$, by two reals 
$a,x$, is presented in which 
$\rL=\HOD 
\;\sneq\;
\rL[a]=\HNT
\;\sneq\;
\rV=\rL[a,x]\,,
$ 
so that $\HNT$ is a model of $\ZFC$ strictly between 
$\HOD$ and the universe.
\end{abstract}

\maketitle

\punk{Introduction}
\las{int}

A set $x$ is  
\rit{nontypical with a cardinal parameter $\ka$}, 
for short $x\in\NT_\ka$,
if it belongs to an 
\OD\ (ordinal definable) set $X$ of cardinality 
$\car X<\ka$.
A set $x$ is  
\rit{hereditarily nontypical 
with a cardinal parameter $\ka$}, 
for short $x\in\HNT_\ka$,
if it itself, all its elements, elements of elements,
and so on, are all nontypical, in other words
the \rit{transitive closure} $\tc(x)$ satisfies
$\tc(x)\sq\NT_\ka$.
These notions Tzouvaras \cite{tz20,tz21} connected 
with some  philosophical and mathematical ideas of 
Bertrand Russell and works of van Lambalgen 
\cite{vL} \etc\ on the concept of randomness. 
They contribute to the ongoing study of 
important classes of sets in the 
set theoretic universe $\rV$ which themselves
satisfy the axioms of set theory, similarly to the  
G\"odel class $\rL$ of all \rit{constructible} sets  
and the class
$\HOD$ of all \rit{hereditarily ordinal definable} sets 
\cite{jechmill}. 

It is clear that $\NT_2=\OD$ and $\HNT_2=\HOD$, thus the 
case $\ka=2$ corresponds to the ordinal definability. 
The classes $\NT_\om$
(elements of finite ordinal definable sets)
and $\HNT_\om$ correspong to  
\rit{algebraically definability} recently 
studied in \cite{FGH,GH,HL}. 
The following classes  
correspond to the next cardinality level $\ka=\omi$:
$$
\NT :=\NT_{\omi}
\;\qand\;
\HNT :=\HNT_{\omi}\,.
$$
Thus $x\in\NT$ iff $x$ belongs to a countable  
\OD\ set, 
and $x\in\HNT$ iff $\tc(x)\sq\NT$.

The class $\HNT$ is transitive and, as shown
in \cite{tz21}, satisfies all axioms of 
$\ZF$ (the axiom of choice \AC\ not included), 
and also satisfies
the relation $\HOD\sq\HNT\sq\rV.$ 
Tzouvaras \cite[2.15]{tz21} asks whether the double 
strict inequality $\HOD\sneq\HNT\sneq\rV$ 
can be realized in an appropriate model of $\ZFC$. 
The following theorem, the main result of this paper, 
answers this question in the affirmative.

\bte
\lam{mt}
Let\/ $\dC=\nse$ be the Cohen forcing for adding 
a generic real\/ $x\in\bn$ to\/ $\rL$. 
There is a forcing notion\/ $\bP\in\rL$, 
which consists of Silver trees, and such that if 
a pair of reals\/ $\ang{a,x}$ is\/ 
\dd{(\bP\ti\dC)}generic over\/ $\rL$ then it is true 
in\/ $\rL[a,x]$ that\vim
$$
\rL=\HOD 
\;\sneq\;
\rL[a]=\HNT
\;\sneq\;
\rV=\rL[a,x]\,.
$$
\ete

Note that $\HNT$ satisfies \ZFC, not 
merely \ZF, in the model $\rL[a,x]$ of the theorem.

\bre
This result is an essential strengthening of 
\cite[Theorem 9.1]{kl63x}. 
Comparably to the latter, the claims that $\rL=\HOD$ 
(instead of simply $a\nin\HOD$) 
and especially $\rL[a]=\HNT$ 
(instead of just $x\nin\HNT$) 
are added here, \poo\ basically the same model, which makes 
the research more accomplished. 

To make the text of this preprint more 
self-contained, we decided to near-copypast 
some definitions and auxiliary results from \cite{kl63x}, 
instead of 
briefly citing them as it 
would be more accustomed in a journal paper. 
\ere

\vyk{
Theorem~\ref{mt} is an essential strengthening of 
\cite[Theorem 9.1]{kl63x}. 
Comparably to the latter, the claims that $\rL=\HOD$ 
(instead of simply $a\nin\HOD$) 
and especially $\rL[a]=\HNT$ 
(instead of just $x\nin\HNT$) 
are added here, \poo\ basically the same model, 
which makes the result more accomplished. 
}

To prove the theorem, we make use of a forcing notion 
$\bP$ introduced in \cite{kl22} in order to define a 
generic real $a\in\dn$ whose \dd{\Eo}equivalence 
class $\eko a$ is a lightface $\ip12$ (hence \OD) set 
of reals with no \OD\ element. 
This property of $\bP$ is responsible for a  
\dd\bP generic real $a$ to belong to $\HNT$, 
and ultimately to $\rL[a]\sq\HNT$, in $\rL[a,x]$. 
This will be based on some results on \sit s and Borel 
functions in Sections \ref{sst},\ref{bnf},\ref{nf}. 
The construction of $\bP$ in $\rL$ is given 
in Sections~\ref{x5},\ref{Uext}. 
The proof that $\rL[a]\sq\HNT$ in $\rL[a,x]$ follows 
in Section~\ref{pa2}. 

The inverse inclusion 
$\HNT\sq\rL[a]$ in $\rL[a,x]$ will be proved 
in Section \ref{pa3} 
on the basis of our earlier result \cite{kl31} 
on countable \OD\ sets in Cohen-generic extensions.

\punk{Perfect trees and \sit{s}}
\las{sst}

Our results will involve forcing notions 
that consist of perfect trees and \sit{s}.
Here we introduce the relevant
terminology from our earlier
works \cite{kl22,kl34,kl30}.

By $\bse$ we denote the set of all \rit{tuples}
(finite sequences) of terms $0,1$,
including the empty tuple $\La$.
The length of a tuple $s$ is denoted by $\lh s$,
and
$2^n=\ens{s\in\bse}{\lh s=n}$ 
(all tuples of length $n$).
A tree $\pu\ne T\sq\bse$ is \rit{perfect},
symbolically $T\in\pet$,
if it has no endpoints and isolated branches.
In this case, the set
$$
[T]=\ens{a\in\dn}{\kaz n\,(a\res n\in T)}
$$
of all \rit{branches} of $T$  
is a perfect set in $\dn.$
Note that $[S]\cap[T]=\pu$ iff $S\cap T$ is finite.\vhm

\bit
\item
If $u\in T\in\pet$, then
a \rit{portion}
(or a \rit{pruned tree}) 
$T\ret u\in\pet$ is defined by 
$T\ret u=\ens{s\in T}{u\su s\lor s\sq u}$.\vhm

\item
A tree $S\sq T$ is \rit{clopen} in $T$ iff it is 
equal to the union of a finite number of portions of $T$. 
This is equivalent to $[S]$ being 
clopen in $[T]$.\vhm
\eit

A tree\/ $T\sq\bse$ is a {\em Silver tree\/},
symbolically\/ $T\in\pes$, if 
there is an infinite
sequence of tuples $u_k=\uu kT\in\bse,$
such that $T$ consists of all tuples of the form
$$
s=u_0\we i_0\we u_1\we i_1\we u_2\we i_2
\we\dots\we u_{n}\we i_n
$$
and their sub-tuples, where $n<\om$ and $i_k=0,1$.
Then the \rit{stem} $\stem T=\uu0T$ is equal 
to the largest tuple $s\in T$ with
$T=T\ret s$, and  
$[T]$ consists of all infinite sequences
$a=u_0\we i_0\we u_1\we i_1\we u_2\we i_2 \we\dots\in\dn,$
where $i_k=0,1$, $\kaz k$.
Put
\pagebreak[0] 
$$
\oin T{n}=\lh{u_0}+1+\lh{u_1}+1+\dots+\lh{u_{n-1}}
+1+\lh{u_n}\,.
$$
In particular, $\oin T{0}=\lh{u_0}$.
Thus 
$\oi T=\ens{\oin Tn}{n<\om}\sq\om$ is the set of all
\rit{splitting levels} of the Silver tree $T$. 

{\ubf Action.} 
Let $\sg\in\bse.$  
If $v\in\bse$ is another tuple of length
$\lh v\ge \lh\sg$, 
then the tuple $v'=\sg\aq v$ of the same length
$\lh {v'}=\lh v$ is defined by 
$v'(i)=v(i)+_2 \sg(i)$
(addition modulo $2$) for all $i<\lh\sg$,
but $v'(i)=v(i)$ whenever $\lh\sg\le i<\lh v$.
If $\lh v<\lh\sg$,
then we just define $\sg\aq v=(\sg\res\lh v)\aq v$.

If $a\in\dn,$ then similarly $a'=\sg\aq a\in\dn,$
$a'(i)=a(i)+_2\sg(i)$ for $i<\lh\sg$,
but $a'(i)=a(i)$ for $i\ge \lh\sg$.
If $T\sq\bse\yt X\sq\dn$, then the sets 
$$
\sg\aq T=\ens{\sg\aq v}{v\in T}
\qand
\sg\aq X=\ens{\sg\aq a}{a\in X}
$$
are \rit{shifts} of the tree $T$ and the set $X$  
accordingly.

\ble
[\cite{kl30}, 3.4]
\lam{LL}
If\/ $n< \om$ and\/ $u,v\in T \cap 2^n,$ then\/
$T \ret u=v \aq u \aq (T\ret v)$.

If\/ $t\in T\in\pes$ and\/ $\sg\in\bse,$ then\/
$\sg\aq T\in\pes$ and\/ $T\ret s\in\pes$.\qed
\ele

\bdf
[\ubf refinements] 
\lam{dr}
Assume that $T,S\in\pes$, $S\sq T$, $n<\om$. 
We define $S\nq n T$
(the tree $S$ \rit{\dd nrefines\/ $T$}) 
if $S\sq T$ and $\oin Tk=\oin Sk$ for all $k<n$. 
This is equivalent to ($S\sq T$ and) 
$\uu kS=\uu kT$ for all $k<n$, of course.
\edf

Then $S\nq 0 T$ is equivalent to $S\sq T$,
and $S\nq{n+1} T$ implies $S\nq n T$ (and $S\sq T$),
but if $n\ge 1$ then $S\nq n T$ is equivalent to
$\oin T{n-1}=\oin S{n-1}$.

\ble
\lam{tadd}
Assume that\/ 
$T,U\in\pes\yt n<\om\yt h>\oin T{n-1} \yt s_0\in 2^h\cap T$,
and\/ $U\sq\req T{s_0}$. 
Then there is a unique tree\/ $S\in\pes$ such that\/
$S\nq n T$ and\/ $\req{S}{s_0}=U.$ 

If in addition\/ $U$ is clopen in\/ $T$ then\/ $S$ 
is clopen in\/ $T$ as well.
\ele
\bpf[sketch]
Define a tree $S$ so that 
$S\cap 2^h=T\cap 2^h$, and if
$t\in T\cap2^h$ then, by Lemma \ref{LL},
$\req {S}t=(t\aq s_0)\aq U$;
then $\req{S}{s_0}=U$. 
To check that $S\in\pes$, we can easily compute the 
tuples $\uu kS$. 
Namely, as $U\sq\req T{s_0}$, we have 
$s_0\sq \uu0U=\roo U$, hence 
$\ell=\lh{(\uu 0U)}\ge h>m=\oin T{n-1}$. 
%
Then $\uu kS=\uu kT$ for all $k<n$, 
$\uu nS=\uu 0U\res[m,\ell)$ 
(thus $\uu nS\in 2^{\ell-m}$), 
and $\uu kS=\uu kU$ for all $k>n$.
%
\epf

\ble
[\cite{kl30}, Lemma 4.4]
\lam{fus}
Let\/
$\ldots \nq 4 T_3\nq 3 T_2\nq 2 T_1\nq 1 T_0$ be 
a sequence of trees in\/ $\pes$.
Then\/ 
$T=\bigcap_nt_n\in\pes$.\qed
\ele
\bpf[sketch]
By definition we have $\uu k{T_n}=\uu k{T_{n+1}}$ 
for all $k\le n$. 
Then one easily computes that $\uu nT=\uu n{T_n}$ for 
all $n$.
\epf

\punk{Reduction of  Borel maps to continuous ones}
\las{bnf}

A classical theorem claims that in Polish
spaces every Borel function is continuous on
a suitable dense $\Gd$ set
(Theorem 8.38 in Kechris \cite{Kdst}). 
It is also known that a Borel map defined on $\dn$ 
is continuous on a suitable \sit. 
The next lemma combines these two results. 
Our interest in functions defined on $\dn\ti\bn$ is 
motivated by further applications to reals in 
generic extensions of the form $\rL[a,x]$, where 
$a\in\dn$ is \dd\bP generic real for some $\bP\sq\pes$ 
while $x\in\bn$ is just Cohen generic.

In the remainder, if $v\in\nse$
(a tuple of natural numbers), then we define 
$\bi v=\ens{x\in\bn}{v\su x}$, 
a {\em Baire interval\/} or {\em portion\/} 
in the Baire space $\bn.$ 

\ble
\lam{su}
Let\/  $T\in\pes$ and\/ 
$f:\dn\ti \bn\to\dn$ be a Borel map.
There is a \sit\/ $S\sq T$ and 
a dense\/ $\Gd$ set\/ $D\sq\bn$
such that\/ $f$ is continuous on\/ $[S]\ti D.$
\ele
\bpf
By the abovementioned classical theorem,
$f$ is already continuous on some dense
$\Gd$ set $Z\sq[T]\ti\bn.$
It remains to define a \sit\/ $S\sq T$ and
a dense $\Gd$ set
$D\sq\bn$ such that $[S]\ti D\sq Z.$
This will be our goal.
 
We have $Z=\bigcap_nZ_n$,
where each $Z_n\sq [T]\ti\bn$ is open dense.

We fix a recursive enumeration
$\om\ti\nse=\ens{\ang{N_k,v_k}}{k<\om}$. 
We will define a sequence of \sit s $S_k$ and 
tuples $w_k\in\nse$ satisfying the following:
\ben
\nenu
\itlb{su1}%
$\ldots \nq 4 S_3\nq 3 S_2\nq 2 s_1\nq 1 S_0=T$,
as in Lemma \ref{fus};

\itlb{su2}%
if $k<\om$ then $S_{k+1}$ is
{\ubf clopen} in $S_k$ (see Section \ref{sst});

\itlb{su3}%
$v_k\sq w_k$ and
$[S_{k+1}]\ti \bi{w_k}\sq Z_{N_k}$, for all $k$.
\een
At step $0$ we already have $S_0=T$.

Assume that $S_k\in\pes$ has already been defined. 
Let $h=\oin{S_k}{k+1}$. 

Consider any tuple $t\in 2^{h}\cap S_k.$
As $Z_{N_k}$ is open dense, there is a tuple 
$u_1\in\nse$ and a \sit\ $A_1\sq\req{S_k}t$,
clopen in $S_k$ (for example, a portion in $S_k$)
such that $v_k\sq u_1$ and
$[A_1]\ti\bi{u_1}\sq Z_{N_k}$.
According to Lemma \ref{tadd}, there exists
a \sit\ $U_1\nq{k+1}S_k$,
clopen in $S_k$ along with $A$,
such that $\req{U_1}t=A_1$, so
$[\req{U_1}t]\ti\bi{u_1}\sq Z_{N_k}$ 
by construction. 

Now take another tuple $t'\in 2^{h}\cap S_k,$
and similarly find 
$u_2\in\nse$ and a \sit\ $A_2\sq\req{U_1}{t'}$,
clopen in $U_1$,
such that $u_1\sq u_2$ and
$[A_2]\ti\bi{u_2}\sq Z_{N_k}$.
Once again there is a \sit\ $U_2\nq{k+1}U_1$,
clopen in $S_k$ and such that 
$[\req{U_2}{t'}]\ti\bi{u_2}\sq Z_{N_k}$.  

We iterate this construction over all tuples 
$t\in 2^{h}\cap S_k,$
\dd{\nq{k+1}}shrinking trees and extending tuples
in $\nse.$
We get a \sit\ $U\nq{k+1}S_k$, clopen in $S_k$, 
and a tuple $w\in\nse,$ that $v_k\sq w$ and
$[U]\ti\bi{w}\sq Z_{N_k}$.
Take $w_k=w\yt S_{k+1}=U$.
This completes the inductive step. 

As a result we get a sequence 
$\ldots \nq 4 S_3\nq 3 S_2\nq 2 S_1 \nq 1 S_0=T$
of \sit{s} $S_k$,
and tuples $w_k\in\nse$ ($k<\om$),
which satisfy \ref{su1},\ref{su2},\ref{su3}.

We put $S=\bigcap_k S_k$; then $S\in\pes$
by \ref{su1} and Lemma \ref{fus}, and $S\sq T$.

If $n<\om$ then let $W_{n}=\ens{w_k}{N_k=n}$.
We claim that $D_{n}=\bigcup_{w\in W_{n}}\bi w$
is an open dense set in $\bn.$
Indeed, let $v\in\nse.$ 
Consider any $k$ such that that $v_k=v$ and $N_k=n$.
By construction, we have $v\sq w_k\in W_{n}$,
as required.
We conclude that the set
$D=\bigcap_{n}D_{n}$ is dense and $\Gd$.

To check  $[S]\ti D\sq Z$,
let $n<\om$; we show that $[S]\ti D\sq Z_n$.
Let $a\in[S]$ and $x\in D$, in particular
$x\in D_n$, so $x\in\bi{w_k}$ for some
$k$ with $N_k=n$.
However, $[S_{k+1}]\ti\bi{w_k}\sq Z_n$ by \ref{su3}, 
and at the same time obviously $a\in[S_{k+1}]$.
We conclude that in fact
$\ang{a,x}\in Z_n$, as required.
\epF{Lemma \ref{su}}

\vyk{
We claim that there exist:\vhm
\ben
\Aenu
\itlb{step4}%
a tuple $w_k\in\nse$ and a \sit\ ,
such that $v_k\sq w_k$ and
$[S_{k+1}]\ti \bi{w_k}\sq Z_{N_k}$.\vhm
\een
}

\bcor
\lam{suC}
Let\/  $T\in\pes$ and 
$f:\dn\to\dn$ be a Borel map.
There is a \sit\/ $S\sq T$ 
such that\/ $f$ is continuous on\/ $[S].$\qed
\ecor

We add the following result that belongs 
to the folklore of the Silver forcing. 
See Corollary 5.4 in \cite{kl34} for a proof.

\ble
\lam{bcL}
Assume that\/  $T\in\pes$ and\/ 
$f:\dn\to\dn$ is a continuous map.
Then there is a \sit\/ $S\sq T$ 
such that\/ $f$ is either a bijection or a constant 
on\/ $[S]$.\qed
\ele

\punk{Normalization of Borel maps}
\las{nf}

\bdf
\lam{nor}
A map $f:\dn\ti\bn\to\dn$ is 
\rit{normalized on\/ $T\in\pes$ for\/ $\bU\sq\pes$}
if there exists a dense\/ $\Gd$ set\/
$X\sq\bn$ such that $f$ is continuous on 
$[T]\ti X$ and$:$\vom
\bit
\item[$-$]
either $(a)$
there are tuples 
$v\in\nse\yt \sg\in\bse$ 
such that\/ $f(a,x)=\sg\aq a$ for all 
$a\in [T]$ and $x\in\bi v\cap X$, where, we remind, 
$\bi v=\ens{x\in\bn}{v\su x}\,;$\vom

\item[$-$]
or $(b)$
$f(a,x)\nin\bigcup_{\sg\in\bse\land S\in \bU}\sg\aq[S]$ 
for all $a\in [T]$ and $x\in X.$\qed
\eit
\eDf

\bte
\lam{tn}
Let\/ $\bU=\ans{T_0,T_1,T_2,\dots}\sq\pes$ and\/
$f:\dn\ti\bn\to\dn $ be a Borel map.
There is a set\/
$\bU'=\ans{S_0,S_1,S_2,\dots}\sq\pes$,   
such that\/ $S_n\sq T_n$ for all\/ $n$ and\/
$f$ is normalized on\/ $S_0$ for\/ $\bU'$.
\ete

\bpf
First of all, according to Lemma \ref{su}, there
is a \sit\ $T'\sq T_0$ and a dense ${\Gd}$ set
$W\sq\bn$ such that $f$ is continuous on $[T']\ti W$.
And since any dense
$\Gd$ set $X\sq\bn$ is homeomorphic to $\bn,$
we can \noo\ assume  
that $W=\bn$ and $T'=T_0$.
Thus, we simply suppose that $f$ is already
\rit{continuous on $[T_0]\ti\bn.$}

Assume that option (a) of the definition of \ref{nor}
does not take place, \ie\vom
\ben
\fenu
\itlb{*}%
if $X\sq\bn$ is dense $\Gd$, and
$v\in\nse\yt \sg\in\bse\yt S\in\pes\yt S\sq T_0$, 
then there are reals
$a\in [S]$ and $x\in\bi v\cap X$ such that 
$f(a,x)\ne\sg\aq a$.\vom
\een
We'll construct \sit s $S_n\sq T_n$ and a dense
$\Gd$ set $X\sq\bn$ satisfying (b) of 
Definition \ref{nor}, that is, in our case,
the relation
$f(a,x)\nin\bigcup_{\sg\in\bse\land n<\om}\sg\aq[S_n]$ 
will be fulfilled 
for all $a\in [S_0]$ and $x\in X.$ 

To maintain the construction, we fix any enumeration
$\om\ti\bse\ti\nse
=\ens{\ang{N_k,\sg_k,v_k}}{k<\om}$. 
Auxiliary \sit s 
$S^n_k$ ($n,k<\om$) and tuples $w_k\in\nse$ ($k<\om$), 
satisfying the following conditions, will be defined. 

\ben
\nenu
\itlb{tn1}%
$\ldots \nq 4 S^n_3\nq 3 S^n_2\nq 2 S^n_1\nq 1 S^n_0=T_n$ 
as in Lemma \ref{fus}, for each $n<\om$;

\itlb{tn2}%
$S^n_{k+1}=S^n_k$ for all $n>0\yt n\ne N_k$;

\itlb{tn3}%
$S^0_{k+1}\nq{k+1}S^0_k$, 
$S^{N_k}_{k+1}\nq{k+1}S^{N_k}_k$,  
$v_k\sq w_k$, and
$f(a,x)\nin\sg_k\aq[S^N_{k+1}]$ 
for all reals $a\in [S^0_{k+1}]$ and $x\in\bi {w_k}$.
\een

At step $0$  of the construction, we put 
$S^n_0=T_n$ for all $n$, by \ref{tn1}.

Assume that $k<\om$ and all \sit s
$S^n_k\yt n<\om$ are already defined. 
We put $S^n_{k+1}=S^n_k$ for all $n>0\yt n\ne N_k$, 
by \ref{tn2}.

To define the trees $S^0_{k+1}$ and $S^{N_k}_{k+1}$, 
we put $h=\oin{S^0_{k}}{k+1}$, 
$m=\oin{S^N_{k}}{k+1}$.\vom

\rit{Case 1\/}: $N_k>0$.
Take any pair of tuples $s\in2^{h}\cap S^0_{k}$, 
$t\in2^{m}\cap S^{N_k}_{k}$
and any reals $a_0\in[\req{S^0_k}s]$ and $x_0\in\bn.$
Consider any real $b_0\in[\req{S^{N_k}_k}t]$ not equal
to $\sg_k\aq f(a_0,x_0)$.
Let's say $b_0(\ell)=i\ne j= (\sg_k\aq f(a_0,x_0))(\ell)$,
where $i,j\le1\yt \ell<\om$.
By the continuity of $f$, there is a tuple 
$u_1\in\nse$ and \sit\ $A\sq\req{S^0_k}s$
such that $v_k\sq u_1\su x_0$, $a_0\in[A]$, and
$(\sg_k\aq f(a,x))(\ell)=j$ for all $x\in\bi{u_1}$
and $a\in [A]$.
It is also clear that
$B=\ens{\tau\in \req{S^{N_k}_k}t}
{\lh \tau\le\ell\lor \tau(\ell)=i}$ 
is a \sit\ containing $b_0$, 
and $b(\ell)=i$ for all $b\in[B]$.
According to Lemma \ref{tadd}, there are  
\sit{s} $U_1\nq{k+1}S^0_k$ and $V_1\nq{k+1}S^{N_k}_k$,
such that $\req{U_1}s=A$ and $\req{V_1}t=B$, hence 
by construction we have
$\sg_k\aq f(a,x)\nin [\req{V_1}t]$ for all
$a\in[\req{U_1}s]$ and $x\in\bi{u_1}$.

Now consider another pair of tuples 
$s\in2^{h}\cap S^0_{k}$, $t\in2^{m}\cap S^{N_k}_{k}$.
We similarly get Silver trees
$U_2\nq{k+1} U_1$ and $V_2\nq{k+1}V_1$, and a tuple
$u_2\in\nse,$ such that $u_1\sq u_2$ and
$\sg_k\aq f(a,x)\nin [\raw{V_2}{t'}]$ for all
$a\in[\req{U_2}{s'}]$ and $x\in\bi{u_2}$.
In this case, we have $\req{V_2}{t}\sq\req{V_1}{t}$
and $\req{U_2}{s}\sq\req{U_1}{s}$, so that what has already
been achieved at the previous step is preserved.

We iterate through all pairs of 
$s\in2^{h}\cap S^0_{k}$, $t\in2^{m}\cap S^{N_k}_{k}$,
\dd{\nq{k+1}}shrinking trees and extending tuples
in $\nse$ at each step.
This results in a pair of Silver trees
$U\nq{k+1}S^0_k\yt V\nq{k+1}S^{N_k}_k$ and a tuple
$w\in\nse$ such that $v_k\sq w$ and
$\sg_k\hspace{0.1ex}\aq f(a,x)\nin [V]$ for all
reals $a\in[U]$ and $x\in\bi{w}$.
Now to fulfill \ref{tn3}, take $w_k=w$,
$S^0_{k+1}=U,$  and $S^{N_k}_{k+1}=V.$ 
Recall that here $N_k>0$.\vom

\rit{Case 2\/}: $N_k=0$.
Here the construction somewhat changes, and
hypothesis \ref{*} will be used. 
We claim that there exist:
\ben
\nenu
\atc
\atc
\atc
\itlb{step0}%
a tuple $w_k\in\nse$ and a Silver tree
$S^0_{k+1}\nq{k+1}S^0_k$
such that $v_k\sq w_k$ and
$f(a,x)\nin\sg_k\aq[S^0_{k+1}]$ 
for all $a\in [S^0_{k+1}]$, $x\in\bi {w_k}$. 
({\ubf Equivalent to \ref{tn3} as $N_k=0$}.) 
\een
Take any pair of tuples
$s,t\in2^{h}\cap S^0_{k}$, where  
$h=\oin{S^0_k}{k+1}$ as above.
Thus 
$S^0_k\ret{t}=t\aq s\aq (S^0_k\ret{s})$,
by Lemma~\ref{LL}.
According to \ref{*}, there are reals $x_0\in\bi{v}$ and
$a_0\in [S^0_k\ret{s}]$ satisfying
$f(a_0,x_0)\ne \sg_k\aq s\aq t\aq a_0$, or 
equivalently, 
$\sg_k\aq f(a_0,x_0)\ne s\aq t\aq a_0$. 

Similarly to Case 1, we have
$(\sg_k\aq f(a_0,x_0))(\ell)=i\ne 
j= (s\aq t\aq a_0)(\ell)$ 
for some $\ell<\om$ and $i,j\le 1$.
By the continuity of $f$, there is a tuple 
$u_1\in\nse$ and a \sit\  
$A\sq S^0_k\ret s$, clopen in $S^0_k$, 
such that $v_k\sq u_1\su x_0$, $a_0\in[A]$, and
$(\sg_k\aq f(a,x))(\ell)=j$ but 
$(s\aq t\aq a)(\ell)=j$  
for all $x\in\bi{u_1}$ and $a\in [A]$.
Lemma \ref{tadd} gives us a 
\sit\ $U_1\nq{k+1}S^0_k$, clopen in $S^0_k$ as well,  
such that $\req{U_1}s=A$ ---
and then $\req{U_1}t=s\aq t\aq A$.
Therefore 
$\sg_k\aq f(a,x)\nin [\req{U_1}t]$ 
holds for all
$a\in[\req{U_1}s]$ and $x\in\bi{u_1}$ 
by construction. 

Having worked out all pairs of tuples 
$s,t\in2^{h}\cap S^0_{k}$,
we obtain a \sit\ $U\nq{k+1}S^0_k$ and a tuple
$w\in\nse,$ such that $v_k\sq w$ and
$\sg_k\aq f(a,x)\nin[U]$ for all
$a\in[U]$ and $x\in\bi{w}$.
Now to fulfill \ref{step0}, take $w_k=w$
and $S^0_{k+1}=U$.

To conclude, we have for each $n$ a sequence
$\ldots \nq 4 S^n_3\nq 3 S^n_2\nq 2 S^n_1\nq 1 S^n_0=T_n$
of \sit{s} $S^n_k$,
along with tuples $w_k\in\nse$ ($k<\om$),
and these sequences satisfy the requirements
\ref{tn1},\ref{tn2},\ref{tn3} 
(equivalent to \ref{step0} in case $N_k=0$).

We put $S_n=\bigcap_k S^n_k$. 
Then $S_n\in\pes$
by Lemma \ref{fus}, and $S_n\sq T_n$.

If $n<\om$ and $\sg\in\bse$ then let
$W_{n\sg}=\ens{w_k}{N_k=n\land \sg_k=\sg}$.  
The set $X_{n\sg}=\bigcup_{w\in W_{n\sg}}\bi w$ 
is then open dense in $\bn.$
Indeed, if $v\in\bn$ then we take $k$
such that $v_k=v\yt N_k=n\yt\sg_k=\sg$;
then $v\sq w_k\in W_{n\sg}$ by construction.
Therefore,
$X=\bigcap_{n<\om\yi \sg\in\bse}X_{n\sg}$  
is a dense $\Gd$ set.
Now to check property (b) of Definition \ref{nor}, 
consider 
any $n<\om\yt\sg\in\bse\yt a\in[S_0]\yt x\in X$;
we claim that $f(a,x)\nin\sg\aq[S_n]$.

By construction, we have $x\in X_{n\sg}$, \ie\
$x\in\bi{w_k}$, where $k\in W_{n\sg}$, so that
$N_k=n\yt \sg_k=\sg$.
Now $f(a,x)\nin\sg\aq[S_n]$ directly
follows from \ref{tn3} for this $k$, since
$S_0\sq S^0_{k+1}$ and $S_n\sq S^n_{k+1}$.
\epF{Theorem~\ref{tn}}

\punk{The  forcing notion for Theorem \ref{mt}}
\las{x5}

Using the standard encoding of Borel sets,
as \eg\ in \cite{sol}
or \cite[\S\,1D]{kl1}, we fix
a coding of Borel functions $f:\dn\to \dn.$
As usual, it includes a 
\rit{\dd{\ip11}set\snos
{The letters $\is{}{}$ and $\ip{}{}$ denote
effective (lightface) projective classes.}
of codes\/} 
$\bc\sq\bn$,
and for each code $r\in\bc$ a certain 
Borel function $F_r:\dn\to \dn$ coded by $r.$
We assume that each Borel function has some
code, and there is a 
$\is11$ relation $\gS(\cdot,\cdot,\cdot)$ and
a $\ip11$ relation $\gP(\cdot,\cdot,\cdot)$
such that for all $r\in\bc$ and $a,b\in\dn$
it holds 
$F_r(a)=b\eqv\gS(r,a,b)\eqv\gP(r,a,b)$.

Similarly, we fix
a coding of Borel functions $f:\dn\ti\bn\to \dn,$
that includes a 
\rit{\dd{\ip11}set of codes\/} 
$\bcd\sq\bn$,
and for each code $r\in\bcd$ a   
Borel function $\fb_r:\dn\ti\bn\to \dn$
coded by $r$, such that each Borel function has some
code, and there is a 
$\is11$ relation $\gsd(\cdot,\cdot,\cdot,\cdot)$ and
a $\ip11$ relation $\gpd(\cdot,\cdot,\cdot,\cdot)$
such that for all $r\in\bc\yt x\in\bn,$ and $a,b\in\dn$
it holds 
$\fb_r(a,x)=b\eqv\gsd(r,a,x,b)\eqv\gpd(r,a,x,b)$.

If $\bU\sq\pes$, then $\clo\bU$ denotes
the set of all trees of the form $\sg\aq(T\ret s)$,
where $\sg\in\bse$ and $s\in T\in\bU$, \ie\ the closure
of $\bU$ \poo\ both shifts and portions.\vhm

{\ubf The following construction is maintained in $\rL$.}
We define a sequence of countable sets
$\bU_\al\sq\pes\yt\al<\omi$ satisfying the following 
conditions \ref{b1}--\ref{b5}.\vhm
\ben
\cenu
\itlb{b1}%
Each $\bU_\al\sq\pes$ is countable,
$\bU_0$ consists of a single
tree $\bse.$\vhm
\een
We then define 
$\bP_\al=\clo{\bU_\al}$,
$\bpl\al=\bigcup_{\xi<\al}\bP_\xi$. 
These
sets are obviously closed with respect
to shifts and portions, that is  
$\clo{\bP_\al}={\bP_\al}$ and
$\clo{\bpl\al}={\bpl\al}$.\vhm
\ben
\cenu
\atc
\itlb{b2}%
For every $T\in\bpl\al$ there is a tree  
$S\in \bU_\al\yt S\sq T$.\vhm
\een
Let $\zfcm$ be the subtheory of the theory \ZFC,
containing all axioms except the power set axiom, and
additionally containing an axiom asserting the existence
of the power set $\pws\om$.
This implies the existence of $\pws X$ for any
countable $X$, the existence of $\omi$ and $\dn$, as well
as the existence of continual sets like $\dn$ or $\pes$.

By $\gM_\al$ we denote the smallest model of 
$\zfcm$ of the form $\rL_\la$
containing the sequence
$\sis{\bU_\xi}{\xi<\al}$, in which $\al$ and all
sets $\bU_\xi\yt\xi<\al$, are countable.\vhm
\ben
\cenu
\atc
\atc
\itlb{b3}%
If a set $D\in\gM_\al\yt D\sq \bpl\al$
is dense in $\bpl\al$, and $U\in\bU_\al$, then
$U\sqf D$,
meaning that there is a finite set $D'\sq D$ such that
$U\sq\bigcup D'$.\vhm

\itlb{b32}%
If a set $D\in\gM_\al\yt D\sq \bpl\al\ti\bpl\al$
is dense in $\bpl\al\ti\bpl\al$, and $U\ne V$ 
belong to $\bU_\al$, then
$U\ti V\sqf D$,
meaning that there is a finite set $D'\sq D$ such that
$[U]\ti[V]\sq\bigcup_{\ang{U',V'}\in D'}[U']\ti[V']$.\vhm
\een
Given that $\clo{\bpl\al}={\bpl\al}$,
this is automatically transferred to all
trees $U\in\bP_\al$ as well.
It follows that $D$ remains
predense in $\bpl\al\cup\bP_\al$.

To formulate the next property, we fix
an enumeration
$$
\pes\ti\bc\ti\bcd=\ens{\ang{T_\xi,b_\xi,c_\xi}}{\xi<\omi}
$$ 
in $\rL$,
which 1) is definable in $\rL_{\omi}$, and
2) each value in $\pes\ti\bc\ti\bcd$ is taken
uncountably many times.\vhm
\ben
\cenu
\atc
\atc
\atc
\atc
\itlb{b4}%
If $T_\al\in\bpl\al$ then there is a tree $S\in\bU_\al$
such that $S\sq T$ and: 
 
(a) $\fb_{b_\al}$
is normalized for $\bU_\al$ on $[S]$
in the sense of Definition \ref{nor}, \ \ and 

(b) $F_{c_\al}$ is 
continuous and either a bijection or a constant on $[S]$.\vhm

\itlb{b5}%
The sequence $\sis{\bU_\al}{\al<\omi}$
is \dd\in definable in $\rL_\omi$.
\een

The construction goes on as follows. 
{\ubf Arguing in $\rL$},
suppose that 
\ben
\fenu
\atc
\itlb{dag}%
$\al<\omi$,
the subsequence $\sis{\bU_\xi}{\xi<\al}$
has been defined and satisfies \ref{b1},\ref{b2} 
below $\al$, and the sets 
$\bP_\xi=\clo{\bU_\xi}$ (for $\xi<\al$), $\bpl\al$,
$\gM_\al$ are defined as above. 
\een
See the proof of the next lemma in Section~\ref{Uext} 
below.

\ble
[$\bU$-extension lemma, in $\rL$]
\lam{susU}%
Under the assumptions of\/ \ref{dag}, 
there is a countable set\/
$\bU_\al\sq\pes$ satisfying\/
\ref{b2}, \ref{b3}, \ref{b32}, \ref{b4}.
\ele

To accomplish the construction, we take $\bU_\al$
to be the smallest, in the sense of the G\"odel 
wellordering of $\rL$,
of those sets that exist by Lemma \ref{susU}.
Since the whole construction is relativized to
$\rL_\omi$, the requirement \ref{b5} is also met.

We put $\bP_\al=\clo{\bU_\al}$ for all $\al<\omi$,
and $\bP=\bigcup_{\al<\omi} \bP_\al$.

The following result, in part related to CCC, 
is a fairly standard
consequence of \ref{b3} and \ref{b32}, see for example
\cite[6.5]{kl22}, \cite[12.4]{kl34}, or
\cite[Lemma 6]{jenmin}; we will skip the proof.

\ble
[in $\rL$]
\lam{cccU}%
The forcing notion\/ $\bP$ belongs to\/ $\rL$, 
satisfies\/ $\bP=\clo{\bP}$ 
and satisfies CCC in\/ $\rL$. 
The product\/ $\bP\ti\bP$ satisfies CCC in\/ $\rL$ 
as well.\qed
\ele

\ble
[in $\rL$]
\lam{norU}%
Assume that\/ $T\in\bP$. 
If \/$g:\dn\to\dn$ is a Borel map 
then there is a tree\/
$S\in\bU_\al$, $S\sq T$, such that\/ $g$
is either a bijection or a constant on\/ $[S]$.

If \/$f:\dn\ti\bn\to\dn $ is a Borel map 
then there is an ordinal\/ $\al<\omi$ and a tree\/
$S\in\bU_\al$, $S\sq T$, such that\/ $f$
is normalized for\/ $\bU_\al$ on\/ $[S]$.
\ele

\bpf
By the choice of the enumeration 
of triples in $\pes\ti\bc\ti\bcd,$
there is an ordinal $\al<\omi$ such that 
$T\in\bpl\al$ and $T=T_\al$, $f=\fb_{b_\al}$, 
$g=F_{b_\al}$.
It remains to refer to \ref{b4}.
\epf

\punk{Proof of the extension lemma}
\las{Uext}

This section is entirely devoted to the 
{\ubf proof of Lemma~\ref{susU}}. 

{\ubf We work in $\rL$ under the assumptions of\/ \ref{dag} 
above.}

We first define a set $\bU=\ens{U_n}{n<\om}$ 
of \sit s $U_n\sq\dn$ satisfying \ref{b2}, \ref{b3} 
\ref{b32}; then further narrowing of the trees 
will be made to also satisfy \ref{b4}. 
This involves a splitting/fusion construction known 
from our earlier papers, see  
\cite[\S\,4]{kl22}, \cite[\S\,9--10]{kl30},
\cite[\S\,10]{kl34}, \cite[\S\,7]{kl49}, 
and to some extent from the proof 
of Theorem~\ref{tn} above. 

We fix enumerations 
$$
\cD=\ens{D(j)}{j<\om}
\qand
\cD_2=\ens{D_2(j)}{j<\om}
$$
of the set $\cD$ of all sets 
$D\in\gM_\al\yt D\sq \bpl\al$
open-dense in $\bpl\al$, and the set $\cD_2$ of all sets 
$D\in\gM_\al\yt D\sq \bpl\al\ti\bpl\al$
open-dense in $\bpl\al\ti\bpl\al$. 
We also fix a bijection $\ba:\om\na\om^4$ which assumes 
each value $\ang{j,j',M,N}\in\om^4$ infinitely many times. 

The construction of the trees $U_n$ is organized 
in the form $U_n=\bigcup_k U^n_k$, where 
the \sit s $U^n_k$ satisfy the following requirements:
\ben
\nenu
\itlb{ux1}%
$\ldots \nq 4 U^n_3\nq 3 U^n_2\nq 2 U^n_1\nq 1 U^n_0$ 
as in Lemma \ref{fus} for each $n<\om$;

\itlb{ux2}%
if $T\in\bpl\al$ then $T=U^n_0$ for some $n$;

\itlb{ux3}%
each $U^n_k$ is a \dd kcollage over $\bpl\al$. 

A \sit\ $T$ is a \dd k\rit{collage} over $\bpl\al$ 
\cite{kl30,kl34} when $T\ret s\in\bpl\al$ for 
each tuple $s\in T\cap 2^h,$ where $h=\oin Tk$. 
Then 0-collages are just trees in $\bpl\al$, 
and every \dd kcollage is a \dd{k+1}collage 
as well since $\clo{\bpl\al}=\bpl\al$.

\itlb{ux4}%
if $k\ge1$, $\ba(k)=\ang{j,j',M,N}$,  
$\mu=\oin{U^M_k}k$, $\nu=\oin{U^N_k}k$ (integers), 
$s\in U^M_k\cap 2^\mu$, 
$t\in U^N_k\cap2^\nu$ (tuples of length resp.\ $\mu,\nu$), 
$M\ne N$, 
then the tree $U^M_k\ret s$ belongs 
to $D(j)$ and the pair $\ang{U^M_k\ret s,U^N_k\ret t}$ 
belongs to $D_2(j')$. 

It follows that $U^M_k\sqf D(j)$ and 
$\ang{U^M_k,U^N_k}\sqf D_2(j')$ in the sense of \ref{b3} 
and \ref{b32} of Section~\ref{x5}.
\een

To begin the inductive construction, we assign 
$U^n_0\in\bpl\al$ so that $\ens{U^n_0}{n<\om}=\bpl\al$, 
to get \ref{ux2}. 
Now let's maintain the step $k\to k+1$. 
Thus suppose that $k<\om$, and all \sit s $U^n_k\yt n<\om$ 
are defined and are \dd kcollages over $\bpl\al$.

Let $\ba(k)=\ang{j,j',M,N}$. 
If $N=M$ then put $U^n_{k+1}=U^n_k$ for all $n$. 

Now assume that $M\ne N.$ 
Put $U^n_{k+1}=U^n_k$ for all $n\nin\ans{M,N}$. 

It takes more effort to define $U^M_{k+1}$ and $U^N_{k+1}$.  
Let $\mu=\oin{U^M_k}{k+1}$, $\nu=\oin{U^N_k}{k+1}$. 
To begin with we put 
$U^M_{k+1}:=U^M_k$ and $U^N_{k+1}:=U^N_k$. 
These \dd{k+1}collages are the initial values for the 
trees $U^M_{k+1}$ and $U^N_{k+1}$, 
to be \dd{\nq{k+1}}shrinked in a finite number of substeps 
(within the step  $k\to k+1$), each corresponding to a 
pair of tuples $s\in U^M_k\cap 2^\mu$ and  
$t\in U^N_k\cap2^\nu$. 

Namely let $s\in U^M_{k+1}\cap 2^\mu$,  
$t\in U^N_{k+1}\cap2^\nu$ be the first such pair. 
The trees $U^M_{k+1}\ret s$, $U^N_{k+1}\ret t$ belong 
to $\bpl\al$ as $U^M_{k+1}$, $U^N_{k+1}$ are 
\dd{k+1}collages over $\bpl\al$. 
Therefore by the open density there exist trees 
$A,B\in D(j)$ such that the pair 
$\ang{U^M_{k+1}\ret s,U^N_{k+1}\ret t}$ belongs to $D_2(j')$ 
and $A\sq U^M_{k+1}\ret s\yt B\sq U^N_{k+1}\ret t$. 
Now Lemma \ref{tadd} gives us \sit s $S\nq{k+1}U^M_k$ 
and $T\nq{k+1}U^N_k$ satisfying $S\ret s\sq A$, 
$T\ret t\sq B$. 
Moreover, by Lemma~\ref{LL}, $S$ and $T$ still are 
\dd{k+1}collages over $\bpl\al$ since $\bpl\al$ 
is closed under shifts by construction. 
To conclude, we have defined 
\dd{k+1}collages 
$S\nq{k+1}U^M_{k+1}$ and $T\nq{k+1}U^N_{k+1}$ 
over $\bpl\al$, satisfying $S\ret s\in D(j)$, 
$T\ret t\in D(j)$, and $\ang{S\ret s,T\ret t}\in D_2(j')$. 
We re-assign the ``new'' $U^M_{k+1}$ and $U^N_{k+1}$ 
to be equal to resp.\ $S,T$. 

Applying this \dd{\nq{k+1}}shrinking procedure 
consecutively for all pairs of tuples 
$s\in U^M_k\cap 2^\mu$ and  
$t\in U^N_k\cap2^\nu$, we eventually 
(after finitely many substeps according to the number 
of all such pairs), we get a pair of \dd{k+1}collages 
$U^M_{k+1}\nq{k+1}U^M_{k}$ and $U^N_{k+1}\nq{k+1}U^N_{k}$ 
over $\bpl\al$, such that for every pair of tuples 
$s\in U^M_k\cap 2^\mu$ and $t\in U^N_k\cap2^\nu$, 
we have $U^M_{k+1}\ret s\in D(j)$  
and 
$\ang{U^M_{k+1}\ret s,U^N_{k+1}\ret t}\in D_2(j')$, 
so conditions 
\ref{ux3} and \ref{ux4} are satisfied. 

Having defined, in $\rL$, a system of \sit s $U^n_k$ 
satisfying \ref{ux1},\ref{ux2},\ref{ux3},\ref{ux4}, we 
then put $U_n=\bigcap_k U^N_k$ for all $n$. 
Those are \sit s by Lemma \ref{fus}. 
The collection $\bU_\al:=\ens{U_n}{n<\om}$ 
satisfies \ref{b2} of Section \ref{x5} by \ref{ux2}. 

To check condition \ref{b3} of Section \ref{x5}, 
let $D\in\gM_\al\yt D\sq \bpl\al$
be dense in $\bpl\al$, and $U\in\bU_\al$. 
We can \noo\ assume that $D$ is open-dense, for if not 
then replace $T$ by  
$D'=\ens{S\in\bpl\al}{\sus T\in D\,(S\sq T)}$. 
Then $D=D(j)$ for some $j$, and $U=U_M$ for some $M$ 
by construction. 
Now consider any index $k$ such that $\ba(k)=\ang{M,N,j,j'}$ 
for $M,j$ as above and any $N,j'$. 
Then we have $U=U_M\sq U^M_k$ by construction, and 
$U^M_k\sqf D$ by \ref{ux4}, thus $U\sqf D$, as required. 

Condition \ref{b32} is verified similarly.  

It remains to somewhat shrink all trees $U_n$ 
to also fulfill \ref{b4}. 
We still work in $\rL$.

Recall that an enumeration
$ 
\pes\ti\bc\ti\bcd=\ens{\ang{T_\xi,b_\xi,c_\xi}}{\xi<\omi}
$,  
parameter-free 
definable in $\rL_{\omi}$, is fixed in Section \ref{x5}. 
We suppose that the tree $T_\al$ belongs to $\bpl\al$.
(If not then we don't worry about \ref{b4}.)  
Consider, according to \ref{b2}, a tree
$U=U_M\in\bU_\al$ satisfying $T\sq T_\al$.
Using Corollary~\ref{suC}, Lemma~\ref{bcL}, 
and Theorem \ref{tn}, we shrink each
tree $U_n\in\bU_\al$ to a tree $U'_n\in\pes\yt U'\sq U$,
so that the function $\fb_{b_\al}$ is normalized
on $U'_M$ for $\bU'=\ens{U'_n}{n<\om}$ 
and $F_{c_\al}$ is 
continuous and either a bijection or a constant on $[U'_M]$.
Take $\bU'$ as the final 
$\bU_\al$ and $T'$ as $U'_M$ to fulfill \ref{b4}.\vom

\qeDD{Lemma~\ref{susU}}

\punk{The model, part I}
\las{x4}

We use the product $\bP \ti\dC$ of the 
forcing notion $\bP$ defined in $\rL$ in Section~\ref{x5} 
and satisfying conditions \ref{b1}--\ref{b5}
as above, and the Cohen forcing,
here in the form of $\dC=\nse$, to prove the 
following more detailed form of Theorem~\ref{mt}. 
The proof of this theorem in the next three 
sections is based on a combination of 
different ideas.

\bte
\lam{nt}
Let a pair of reals\/ $\ang{a_0,x_0}$ be\/
\dd{\bP\ti\dC}generic over \/ $\rL$.
Then 
\ben
\Renu
\itlb{nt1}%
$a_0$ is not $\OD$, and moreover, 
$\HOD=\rL$  in\/ $\rL[a_0,x_0]\,;$

\itlb{nt2}%
$a_0$ belongs to $\HNT$, and moreover, 
$\rL[a_0]\sq\HNT$  in\/ $\rL[a_0,x_0]\,;$

\itlb{nt3}%
$x_0$ does not belong to $\HNT$, and moreover, 
$\HNT\sq\rL[a_0]$  in\/ $\rL[a_0,x_0]\,.$
\een
\ete

We prove Claim~\ref{nt1} of the theorem 
in this section. 
The proof is based on several lemmas. 
According to the next lemma, 
it suffices to prove that $\HOD=\rL$ in $\rL[a_0]$. 

\ble
\lam{I0}
$(\HOD)^{\rL[a_0,x_0]}\sq(\HOD)^{\rL[a_0]}$.
\ele
\bpf
By the forcing product theorem, $x_0$ is a 
Cohen generic real over $\rl[a_0]$. 
It follows by a standard argument based on the full 
homogeneity of the Cohen forcing $\dC$ that 
if $H\sq\Ord$ is $\OD$ in $\rL[a_0,x_0]$ then 
$H\in \rL[a_0]$ and $H$ is $\OD$ in $\rL[a_0]$. 

Now prove the implication  
$Y\in (\HOD)^{\rL[a_0,x_0]}\imp 
Y\in \rL\land Y\in (\HOD)^{\rL[a_0]}$ 
by induction on
the set-theoretic rank $\rk x$ of $x\in\rL[a_0,x_0]$.
Since each set consists only of sets
of strictly lower rank, it is sufficient to check that
if a set $H\in\rL[a_0,x_0]$ satisfies
$H\sq(\HOD)^{\rL[a_0]}$ and $H\in\HOD$ in $\rL[a_0,x_0]$ 
then $H\in\rL[a_0$ and $H\in (\OD)^{\rL[a_0]}$.
Here we can assume that in fact $H\sq\Ord$,
since $\HOD$ allows an \OD\ wellordering and hence an 
\OD\ bijection onto $\Ord$. 
But in this case $H\in \rL[a_0]$ and 
$H$ is $\OD$ in $\rL[a_0]$ by the above, as required. 
\epf

\ble
[Lemma 7.5 in \cite{kl22}]
\lam{I1}
$a_0$ is not\/ $\OD$ in\/ $\rL[a_0]$.
\ele
\bpf
Suppose towards the contrary that $a_0$ is $\OD$ in $\rL[a_0]$.
But $a_0$ is a \dd\bP generic real over $\rL$, so the 
contrary assumption is forced. 
In other words, there is a tree $T\in\bP$ with $a_0\in[T]$ 
and a formula $\vt(x)$ with ordinal parameters, such that 
if $a\in[T]$ is \dd\bP generic over $\rL$ then $a$ is the 
only real in $\rL[a]$ satisfying $\vt(a)$. 
Let $s =\roo(T )$. 
Then both $s\we 0$ and $s\we 1$ belong to $T$, and 
either $s\we 0\su a_0$ or $s\we 1\su a_0$. 
Let, say, $s\we 0\su a_0$. 
Let $n = \lh(s)$ and $\sg = 0^n \we 1$, 
so that all three strings $s\we 0$, $s\we 1$, $\sg$ 
belong to $2^{n+1},$ and $s\we 0=\sg\aq (s\we 1)$. 
As the forcing $\bP$ is invariant under the action of $\sg$, 
the real $a_1=\sg\aq a_0$ is $\bP$-generic over $\rL$, 
and $\sg\aq T=T$.  
We conclude that it is true in $\rL[a_1]=\rL[a_0]$ 
that $a_1$ is still the only real in $\rL[a_1]$ 
satisfying $\vt(a_1)$. 
However obviously $a_1\ne a_0$!
\epf

\ble
\lam{I2}
If\/ $b\in\rL[a_0]\bez\rL$ is a real then\/ $b$ 
is not\/ $\OD$ in\/ $\rL[a_0]$.
\ele
\bpf
It follows from Lemma \ref{cccU} 
(and the countability of $\dC$) 
that the forcing $\bP\ti\dC$ preserves cardinals.
We conclude that that $b=g(a_0)$ for some Borel
function $g=F_r:\dn\to\dn$ with a code $r\in\bc\cap\rL$.
Now by Lemma~\ref{norU} there is a tree $S\in\bP$ such 
that $a_0\in[S]$ and $h=g\res[S]$ is a bijection of a 
constant. 
If $h$ is a bijection then $b\nin\OD$ in $\rL[a_0]$ 
since otherwise 
$a_0=h\obr(b)\in\OD$, contrary to Lemma~\ref{I1}. 
If $h$ is a constant, so that there is a real 
$b_0\in\rL\cap\dn$ such that $h(a)=b_0$ for all $a\in[S]$, 
then $b=h(a_0)=c\in\rL$, contrary to the choice of $b$.
\epf

\ble
\lam{I3}
If\/ $X\sq\Ord\yt X\in\rL[a_0]\bez\rL$, then\/ $X$ 
is not\/ $\OD$ in\/ $\rL[a_0]$.
\ele
\bpf
Suppose to the contrary that $X\sq\Ord$, 
$X\in\rL[a_0]\bez\rL$, and $X$ is $\OD$ in\/ $\rL[a_0]$. 
Let $t$ be a \dd\bP name for $X$. 
Then a condition $T_0\in\bP$ (a \sit) \dd\bP forces 
$$
t\in\rL[a_0]\bez\rL\;\land\; t\in\OD  
$$  
over $\rL$. 
Say that $t$ \rit{splits} conditions $S,T\in\bP$ if 
there is an ordinal $\ga$ suct that 
$S$ forces $\ga\in t$ but $T$ forces $\ga\nin t$ or 
vice versa; let $\ga_{ST}$ be the least such an 
ordinal $\ga$. 

We claim that the set 
$$
D=\ens{\ang{S,T}}
{S,T\in\bP\land S\cup T\sq T_0\land 
t\text{ splits }S,T}\in\rL
$$
is dense in $\bP\ti\bP$ above $\ang{T_0,T_0}$. 
Indeed let $S,T\in\bP$ be subtrees of $T_0$. 
If $t$ splits no stronger pair of trees $S'\sq S$, 
$T'\sq T$ in $\bP$ then easily both $S$ and $T$ 
decide $\ga\in t$ for every ordinal $\ga$, a 
contradiction with the choice of $T_0$. 
Thus $D$ is indeed dense. 

Let, in $\rL$, $A\sq D$ be a maximal antichain; 
$A$ is countable in $\rL$ by Lemma~\ref{cccU}, 
and hence the set 
$W=\ens{\ga_{ST}}{\ang{S,T}\in A}\in\rL$ is countable 
in $\rL$. 
We claim that 
\ben
\fenu
\atc
\atc
\itlb{**}%
the intersection $b=X\cap W$ does not belong to $\rL$. 
\een
Indeed otherwise there is a tree $T_1\in\bP\yt T_1\sq T_0$, 
which \dd\bP forces that $t\cap W=b$. 
(The sets $W,b\in\rL$ are identified with their names.) 
 
By the countability of $A,W$ there is an ordinal 
$\al<\omil$ such that $A\sq \bpl\al\ti\bpl\al$, 
$T_1\in\bpl\al$, and $W\sq\al$. 
We can \noo\ assume that $A\in \gM_\al$, for if not then 
further increase $\al$ below $\omil$ accordingly. 
Let $u=\roo{T_1}$. 
The trees $T_{10}=T_1\ret{u\we0}$ and 
$T_{11}=T_1\ret{u\we1}$ belong to $\bpl\al$ along with 
$T_1$, and hence there are trees $U\yt V\in\bU_\al$ 
with $U\sq T_{10}$ and $V\sq T_{11}$. 
Clearly $U\ne V$, so that we have 
$[U]\ti [V]\sq\bigcup_{\ang{U',V'}\in A'}[U']\ti[V']$ 
for a finite set $A'\sq A$ by \ref{b32} of Section \ref{x5}. 
Now take reals $a'\in [U]$ and $a''\in[V]$ both 
\dd\bP generic over $\rL$. 
Then there is a pair of trees $\ang{U',V'}\in A'$ such 
that $a'\in [U']$ and $a''\in[V']$. 
The interpretations $X'=t[a']$ and $X''=t[a'']$ are then 
different on the ordinal $\ga=\ga_{U'U''}\in W$ since 
$A'\sq A\sq D$. 
Thus the restricted sets $b'=X'\res W$ and $b''=X''\res W$ 
differ from each other. 
In particular at least one of $b',b''$ is not equal to $b.$ 
But $a',a''\in[T_1]$ by construction, hence this contradicts 
the choice of $T_1$ and completes the proof of \ref{**}. 

Recall that $b\sq W,$ and $W\in\rL$ is countable in $\rL$. 
It follows that $b$ can be considered as a real, 
so we conclude that $b$ is not \OD\ in $\rL[a_0]$ 
by Lemma~\ref{I2} and \ref{**}. 

However $b=X\cap W,$ where $X$ is \OD\ and $W\in\rL$, 
hence $W$ is \OD\ in $\rL[a_0]$ and 
$b$  is \OD\ in $\rL[a_0]$. 
The contradiction obtained ends the proof of the lemma.
\epF{Lemma}

Now Theorem~\ref{nt}\ref{nt1} 
immediately follows from 
Lemma~\ref{I0} and Lemma~\ref{I3}.\vtm

\qeDD{Claim~\ref{nt2} of Theorem~\ref{nt}}

\punk{The model, part II}
\las{pa2}

Here we establish Claim~\ref{nt2} of Theorem~\ref{nt}. 
To prove $\rL[a_0]\sq\HNT$ it suffices to show  
that $a_0$ itself belongs to $\HNT$, and then make use 
of the fact that by G\"odel every set $z\in\rL[a_0]$ 
has the form $x=F(a_0)$, where $F$ is an \OD\ function. 

Further, to prove $a_0\in\HNT$
it suffices to check that the 
\dd\Eo equivalence class\snos
{Recall that the equivalence relation 
$\Eo$
is defined on $\dn$ so that $a\Eo b$ iff the set
$a\sd b=\ens{k}{a(k)\ne b(k)}$ is finite. 
Equivalently, $a\Eo b$ iff $a=\sg\aq b$ for some 
tuple $\sg\in\bse$. 
Then 
$\eko a=\ens{b\in\dn}{a\Eo b}=\ens{\sg\aq a}{\sg\in\bse}$ 
is the \dd\Eo\rit{equivalence class} of $a$.}
$\eko{a_0}=\ens{b\in\dn}{a_0\Eo b}$  
(which is a countable set) of our generic 
real $a_0$ is an \OD\ set in $\rL[a_0,x_0]$.
According to \ref{b5}, it suffices to establish
the equality
$$
\textstyle
\eko{a_0}=\bigcap_{\xi<\omi}\bigcup_{T\in\bP_\xi}[T]\,.
\eqno(*)
$$
Note that every set $\bP_\xi$ is pre-dense
in $\bP$; this follows from \ref{b3} and \ref{b4}, see,
for example, Lemma 6.3 in \cite{kl22}.
This immediately implies $a_0\in \bigcup_{T\in\bP_\xi}[T]$
for each $\xi$.
Yet all sets $\bP_\xi$ are invariant \poo\ 
shifts by construction.
Thus we have $\sq$ in (*).

To prove the inverse inclusion, assume
that a real $b\in\dn$ belongs to the right-hand side 
of (*) in $\rL[a_0,x_0]$.
It follows from Lemma \ref{cccU} 
(and the countability of $\dC$) 
that the forcing $\bP\ti\dC$ preserves cardinals.
We conclude that that $b=g(a_0,x_0)$ for some Borel
function $g=F_q:\dn\ti\bn\to\dn$ with a code $q\in\bc\cap\rL$.

{\ubf Assume to the contrary that} $b=g(a_0,x_0)\nin\eko{a_0}$.

Since $x_0\in\bn$ is a \dd\dC generic
real over $\rl[a_0]$ by the forcing product theorem,
this assumption is forced, so that there is a 
tuple $u\in\dC=\nse$ such that
$$
\textstyle
f(a_0,x)\in\bigcap_{\xi<\omi}\bigcup_{T\in\bP_\xi}[T]
\bez \eko{a_0}\,,
$$
whenever a real $x\in\bi u$ is
\dd\dC generic over $\rL[a_0]$.
(Recall that $\bi u=\ens{y\in\bn}{u\su y}$.)
Let $H$ be the canonical homomorphism of 
$\bn$ onto $\bi u$. 
We put $f(a,x)=g(a,H(x))$ for $a\in\dn\yt x\in\bn.$
Then $H$ preserves the \dd\dC genericity,
and hence
$$
\textstyle
f(a_0,x)\in\bigcap_{\xi<\omi}\bigcup_{T\in\bP_\xi}[T]
\bez \eko{a}\,,
\eqno(**)
$$
whenever $x\in\bn$ is 
\dd\dC generic over $\rL[a_0]$.
Note that $f$ also has a Borel code 
$r\in\bc$ in $\rL$, so that $f=F_{r}$.

It follows from Lemma \ref{norU} that there is an ordinal 
$\al<\omi$ and a tree
$S\in\bU_\al$, on which $f$ is normalized for $\bU_\al$, 
and which satisfies $a_0\in [S]$.
Normalization means that, in $\rL$,
there is a dense $\Gd$ set $X\sq\bn$
satisfying one of the two options of 
Definition~\ref{nor}.
Consider a real $z\in\bn\cap\rL$
(a \rit{$\Gd$-code} for $X$ in $\rl$) 
such that 
$X=X_z=\bigcap_k\bigcup_{z(2^k\cdot 3^j)=1}\bi{w_j}$,
where $\bse=\ens{w_j}{j<\om}$ is a 
fixed recursive enumeration of tuples.\vom

\rit{Case 1}:
there are tuples
$v\in\nse\yt \sg\in\bse,$ 
such that $f(a,x)=\sg\aq a$ for all 
points $a\in [S]$ and $x\in\bi v\cap X$.
In other words, it is true in $\rL$ that 
$$
\kaz a\in[S]\,\kaz x\in\bi v\cap X_z\:
(f(a,x)=\sg\aq a)\,.
$$
But this formula is absolute by Shoenfield, so it is also true
in $\rL[a_0,x_0]$.
Take $a=a_0$ (recall: $a_0\in[S]$)
and any real $x\in\bi v$,
\dd\dC generic over $\rL[a_0]$.
Then $x\in X_z$, because $X_z$ is a dense $\Gd$
with a code even from $\rL$.
Thus $f(a_0,x)=\sg\aq a_0\in\eko{a_0}$,
which contradicts (**).\vom

\rit{Case 2}:
$f(a,x)\nin\bigcup_{\sg\in\bse\land U\in \bU_\al}
\sg\aq[U]$ 
for all $a\in [S]$ and $x\in X.$
By the definition of $\bP_\al$, this implies 
$f(a,x)\nin\bigcup_{T\in \bP_\al}[T]$
for all $a\in [S]$ and $x\in X,$
and this again contradicts (**) for $a=a_0$.
\vom

The resulting contradiction in both cases refutes
the contrary assumption above and completes the 
proof.\vom


\qeDD{Claim~\ref{nt2} of Theorem~\ref{nt}}

\punk{The model, part III}
\las{pa3}

Here we prove  Claim~\ref{nt3} of Theorem~\ref{nt}. 
We make use of the following result here. 

\ble
\lam{le3}
Let\/ $x\in\bn$ be Cohen-generic 
over a set universe\/ $\rV.$ 
Then it holds in\/ $\rV[x]$ that if\/
$Z\sq\dn$ is a countable \OD\ set then\/
$Z\in\rV.$ 
More generally if\/ $q\in\dn\cap\rV$  
then it holds in\/ $\rV[x]$ that if\/
$Z\sq\dn$ is a countable\/ $\OD(q)$ set then\/
$Z\in\rV.$ 
\qed
\ele
\bpf[sketch]
The pure \OD\ case is Theorem 1.1 in 
\cite{kl31}.\snos
{See our papers \cite{kl31,kl39,kl40} for more on 
countable and Borel \OD\ sets in Cohen and some 
other generic extensions.} 
The proof of the general case does not differ, 
$q$ is present in the flow of arguments as a passive 
parameter. 
\epf

This result admits the following extension 
for the case  $\rV=\rL$. 
Here $\OD(a)$ naturally means sets definable 
by a formula containing $a_0$ and ordinals as parameters

\bcor
\lam{c3}
Assume that\/ $a\in\dn$ and\/ $x\in\bn$ is Cohen-generic 
over\/ $\rL[a]$.
Then it holds in\/ $\rL[a,x]$ that
if $X\in\rL[a]$ and\/ $A\sq2^X$ is a
countable $\OD(a)$ set then $A\sq\rL$.
\ecor
\bpf
As the Cohen forcing is countable,
there is a set $Y\sq X\yt Y\in\rL[a]$, 
countable in $\rL[a]$ and such that if $f\ne g$
belong to $2^X$ then $f(x)\ne g(x)$ for some
$x\in Y$.
Then $Y$ is countable and $\OD(a)$ in $\rL[a,x]$,
so the \rit{projection}
$B=\ens{f\res Y}{f\in A}$ of the set $A$ will also
be countable and $\OD(a)$ in $\rL[a,x]$.
We have $B\in\rL[a]$ by Lemma~\ref{le3}.
(The set $Y$ here can be identified with $\om$.)
Hence, each $f\in B$ is $\OD(a)$ in $\rL[a,x]$.
However, if $f\in A$ and $w=f\res Y$,
then by the choice of $Y$ it holds in
$\rL[a,x]$ that 
$f$ is the only element in $A$ satisfying 
$f\res Y=w$.
Therefore $f\in\OD(a)$ in $\rL[a,x]$. 
We conclude that $f\in\rL[a]$.
\epf

\bpf[Claim~\ref{nt3} of Theorem~\ref{nt}]
We prove an even stronger claim 
$$
{x\in\HNT(a_0)}\imp{x\in\rL[a_0]}
$$ 
in $\rL[a_0,x_0]$ by induction on
the set-theoretic rank $\rk x$ of 
sets $x\in\rL[a_0,x_0]$. 
Here $\HNT(a_0)$ naturally means all sets hereditarily 
$\NT(a_0)$, the latter meals all elements of countable 
sets in $\OD(a_0)$.

Since each set consists only of sets
of strictly lower rank, it is sufficient to check that
if a set $H\in\rL[a_0,x_0]$ satisfies
$H\sq\rL[a_0]$ and $H\in\HNT(a_0)$ in $\rL[a_0,x_0]$ 
then $H\in\rL[a_0]$.
Here we can assume that in fact $H\sq\Ord$,
since $\rL[a_0]$ allows an $\OD(a_0)$ wellordering.
Thus, let $H\sq\la\in\Ord$.
Additionally, since $H\in\HNT(a_0)$,
we have, in $\rL[a_0,x_0]$, a countable $\OD(a_0)$ set
$A\sq\pws\la$ containing $H$.
However, $A\in\rL[a_0]$ by Corollary \ref{c3}.
This implies $H\in\rL[a_0]$ as required. 
\epf

\qeDD{Claim~\ref{nt3} and Theorem~\ref{nt} as a whole}\vom

\qeDD{Theorem~\ref{mt}}

\punk{Comments and questions}
\las{x7}

\hspace*{0pt}\indent 
1. 
Recall that if $x$ is a Cohen real over $\rL$ 
then $\HNT=\rL$ in $\rL[x]$ by Lemma~\ref{le3}. 

\bqe
\lam{QQ}
Is it true in generic extensions of $\rL$ by a single 
Cohen generic real that a countable \OD\ set of 
any kind necessarily consists 
only of \OD\ elements?
\eqe

We cannot solve this even for \rit{finite} \OD\ sets. 

By the way it is not that obvious to expect the \rit{positive} 
answer. 
Indeed, the problem solves in the \rit{negative} for 
Sacks and some other generic extensions even for \rit{pairs}, 
see \cite{kl54,kl62}. 
For instance, if $x$ is a Sacks-generic real over $\rL$ then 
it is true in $\rL[x]$ that there is an $\OD$ unordered pair 
$\ans{X,Y}$ of sets of reals $X,Y\sq\pws\dn$ such that $X,Y$ 
themselves are non-\OD\ sets. 
See \cite{kl54} for a proof of this rather surprising 
result originally by Solovay.

2. 
See Fuchs \cite{fu} (unpublished) for 
some other research lines 
related to Russell-nontypical sets with various cardinal 
parameters.


\bibliographystyle{amsplain}

\small

\bibliography{64,64kl}

\providecommand{\bysame}{\leavevmode\hbox to3em{\hrulefill}\thinspace}
\providecommand{\MR}{\relax\ifhmode\unskip\space\fi MR }
\providecommand{\MRhref}[2]{%
  \href{http://www.ams.org/mathscinet-getitem?mr=#1}{#2}
}
\providecommand{\href}[2]{#2}
\begin{thebibliography}{10}

\bibitem{kl54}
Ali Enayat and Vladimir Kanovei, \emph{{An unpublished theorem of Solovay on OD
  partitions of reals into two non-OD parts, revisited}}, {J. Math. Log.}
  \textbf{21} (2021), no.~3, 1--22, Article No. 2150014.

\bibitem{kl62}
Ali {Enayat}, Vladimir {Kanovei}, and Vassily {Lyubetsky}, \emph{{On
  effectively indiscernible projective sets and the Leibniz-Mycielski axiom}},
  {Mathematics} \textbf{9} (2021), no.~14, 1--19 (English), Article No 1670.

\bibitem{fu}
Gunter {Fuchs}, \emph{Blurry definability}, Preprint, 2021.

\bibitem{FGH}
Gunter {Fuchs}, Victoria {Gitman}, and Joel~David {Hamkins},
  \emph{{Ehrenfeucht's lemma in set theory}}, {Notre Dame J. Formal Logic}
  \textbf{59} (2018), no.~3, 355--370 (English).

\bibitem{GH}
Marcia~J. {Groszek} and Joel~David {Hamkins}, \emph{{The implicitly
  constructible universe}}, {J. Symb. Log.} \textbf{84} (2019), no.~4,
  1403--1421 (English).

\bibitem{HL}
Joel~David {Hamkins} and Cole {Leahy}, \emph{{Algebraicity and implicit
  definability in set theory}}, {Notre Dame J. Formal Logic} \textbf{57}
  (2016), no.~3, 431--439 (English).

\bibitem{jechmill}
Thomas {Jech}, \emph{{Set theory}}, {The third millennium revised and expanded}
  ed., Springer-Verlag, Berlin-Heidelberg-New York, 2003 (English), Pages xiii
  + 769.

\bibitem{jenmin}
Ronald {Jensen}, \emph{{Definable sets of minimal degree}}, {Math. Logic Found.
  Set Theory, Proc. Int. Colloqu., Jerusalem 1968} (Yehoshua Bar-Hillel, ed.),
  Studies in logic and the foundations of mathematics, vol.~59, North-Holland,
  Amsterdam-London, 1970, pp.~122--128.

\bibitem{kl1}
Vladimir {Kanovei} and Vassily {Lyubetsky}, \emph{{On some classical problems
  in descriptive set theory}}, Russ. Math. Surv. \textbf{58} (2003), no.~5,
  839--927 (Russian, English).

\bibitem{kl22}
\bysame, \emph{{A definable $E_0$ class containing no definable elements}},
  {Arch. Math. Logic} \textbf{54} (2015), no.~5-6, 711--723 (English).

\bibitem{kl31}
\bysame, \emph{{Countable OD sets of reals belong to the ground model}}, {Arch.
  Math. Logic} \textbf{57} (2018), no.~3-4, 285--298 (English).

\bibitem{kl34}
\bysame, \emph{{Definable $\mathsf{E}_0$ classes at arbitrary projective
  levels}}, {Ann. Pure Appl. Logic} \textbf{169} (2018), no.~9, 851--871
  (English).

\bibitem{kl30}
\bysame, \emph{{Non-uniformizable sets of second projective level with
  countable cross-sections in the form of Vitali classes}}, Izvestiya:
  Mathematics \textbf{82} (2018), no.~1, 61--90.

\bibitem{kl40}
\bysame, \emph{{Borel OD sets of reals are OD-Borel in some simple models}},
  {Proc. Am. Math. Soc.} \textbf{147} (2019), no.~3, 1277--1282 (English).

\bibitem{kl39}
\bysame, \emph{{Definable elements of definable Borel sets}}, {Math. Notes}
  \textbf{105} (2019), no.~5, 684--693 (English).

\bibitem{kl49}
\bysame, \emph{{Models of set theory in which separation theorem fails}},
  Izvestiya: Mathematics \textbf{85} (2021), no.~6, to appear.

\bibitem{kl63x}
\bysame, \emph{{On Russell typicality in Set Theory}}, arXiv e-prints (2021),
  arXiv:2111.07654.

\bibitem{Kdst}
Alexander~S. Kechris, \emph{Classical descriptive set theory}, Springer-Verlag,
  New York, 1995. \MR{96e:03057}

\bibitem{vL}
{Michiel van} {Lambalgen}, \emph{{The axiomatization of randomness}}, {J. Symb.
  Log.} \textbf{55} (1990), no.~3, 1143--1167 (English).

\bibitem{sol}
{R}obert~M. Solovay, \emph{A model of set-theory in which every set of reals is
  {L}ebesgue measurable}, Ann. Math. (2) \textbf{92} (1970), 1--56.

\bibitem{tz21}
Athanassios Tzouvaras, \emph{Typicality \'a la {R}ussell in set theory},
  ResearchGate Preprint, May 2021, 14 pp.,
  \href{https://www.researchgate.net/publication/351358980_Typicality_a_la_Russell_in_set_theory}{ResearchGate
  Link}.

\bibitem{tz20}
\bysame, \emph{Russell's typicality as another randomness notion}, Mathematical
  Logic Quarterly \textbf{66} (2020), no.~3, 355--365.

\end{thebibliography}

\end{document}